\documentclass[submission]{dmtcs}
\usepackage[latin1]{inputenc}
\usepackage{graphicx}
\usepackage{amssymb,amsmath}
\usepackage[T1]{fontenc}

\newtheorem{theo}{Theorem\ }[section]
\newtheorem{prop}[theo]{Proposition\ }

\newtheorem{definitio}[theo]{Definition\ }
\newenvironment{defi}[1]{\begin{definitio} {\bf #1} \rm }{\end{definitio}}
\newtheorem{algorith}[theo]{Algorithm\ }

\newtheorem{exemple}[theo]{Example\ }
\newenvironment{exam}[1]{\begin{exemple}{\bf #1} \rm }{\end{exemple}}
\newtheorem{exercice}[theo]{Exercise\ }

\newtheorem{notations}[theo]{Notation\ }

\newtheorem{paragr}[theo]{\hspace{-.8ex}}

\newcommand{\End}{\mathop{\sf End}}

\newcommand{\rt}[1]{{\upharpoonright}_{#1}}

\newcommand{\supp}{\mathop{\sf supp}}
\newcommand{\oc}[1]{{\stackrel{\circ}{#1}}}
\newcommand{\ol}[1]{{\overline{#1}}}
\newcommand{\bom}[1]{{\mbox{\boldmath $ #1 $}}}
\newcommand{\card}{\mathop{\rm card}}

\newcommand{\Rk}{{\par\noindent{\bf Remark:~}}}
\newcommand{\bra}{\langle\,}
\newcommand{\ket}{\,\rangle}
\newcommand{\dd}{\mathop{\rm d}}
\newcommand{\cA}{\mathcal{A}}
\newcommand{\cF}{\mathcal{F}}
\newcommand{\cG}{\mathcal{G}}
\newcommand{\cS}{\mathcal{S}}
\newcommand{\cX}{\mathcal{X}}
\newcommand{\be}{{\bf e}}
\newcommand{\bu}{{\bf u}}
\newcommand{\bv}{{\bf v}}
\newcommand{\BbA}{\mathbb{A}}
\newcommand{\BbE}{\mathbb{E}}
\newcommand{\BbN}{\mathbb{N}}
\newcommand{\BbP}{\mathbb{P}}
\newcommand{\BbR}{\mathbb{R}}
\newcommand{\BbT}{\mathbb{T}}
\newcommand{\BbV}{\mathbb{V}}
\newcommand{\BbX}{\mathbb{X}}
\newcommand{\BbZ}{\mathbb{Z}}

\newcommand{\lgr}{\underline{\mathsf{gr}}(\BbV)}
\newcommand{\ugr}{\ol{\mathsf{gr}}(\BbV)}
\newcommand{\grr}{\mathsf{gr}(\BbV)}
\newcommand{\brr}{\mathsf{br}(\BbV)}
\newcommand{\GLdR}{\mathsf{GL}(d,\BbR)}
\newcommand{\0}{{\bf 0}}
\newcommand{\StrPosInt}{\BbN  }
\newcommand{\PosInt}{\BbZ_+ }
\newcommand{\eqlaw}{\stackrel{\text{\tiny law}}{=}}

\author{Mikhail Menshikov\addressmark{1}{\ }
  \and Dimitri Petritis\addressmark{2}{\ }
  \and Serguei Popov\addressmark{3}\thanks{Partially supported by 
CNPq (302981/02-0)}}
\title[Bindweeds]{Bindweeds or 
random walks 
in random environments on multiplexed trees and their asympotics}
\address{\addressmark{1}
Department of Mathematical Sciences, 
University of Durham,
South Road, Durham DH1 3LE,
United Kingdom\\
 Mikhail.Menshikov@durham.ac.uk\\
\addressmark{2} 
Institut de Recherche
Math\'ematique, Universit\'e de Rennes I and CNRS UMR 6625, 
35042 Rennes Cedex, France\\
Dimitri.Petritis@univ-rennes1.fr\\
\addressmark{3}
Instituto de Matem\'atica e Estat\'\i stica,
Universidade de S\~ao Paulo,
Rua do Mat\~ao 1010, CEP 05508-090,
S\~ao Paulo SP, Brasil\\
popov@ime.usp.br}
\revision{0}
\received{\today}
\revised{}
\accepted{}
\begin{document}
\maketitle
\begin{abstract}
  We report on the asymptotic behaviour of
a new model of random walk, we term 
the bindweed model, evolving in a random environment
on an infinite multiplexed tree.
The term \textit{multiplexed} means that the  model can be
viewed as a nearest neighbours random walk on a tree whose
vertices carry an internal degree of freedom from
the finite set $\{1,\ldots,d\}$, for some integer $d$.
The  consequence of the internal degree of freedom is
an enhancement of the tree graph structure induced 
by the replacement of ordinary edges by multi-edges, indexed
by the set $\{1,\ldots,d\}\times\{1,\ldots,d\}$. This indexing
conveys the information on the internal degree of freedom
of the vertices contiguous to each edge.
The term \textit{random environment} means that the jumping rates
for the random walk are a family of edge-indexed 
random variables, independent
of the natural filtration generated by the
random variables entering in the definition of the
random walk; their joint distribution depends on the index of each component
of the multi-edges. We study the large time asymptotic behaviour
of this random walk and classify it with respect
to positive recurrence or transience in terms
of a specific parameter of the probability distribution
of the jump rates.
This classifying parameter is shown to coincide
with the critical value of a 
matrix-valued multiplicative cascade on the ordinary tree (\textit{i.e.}\
the one without internal degrees of freedom attached to the vertices)
having the same vertex set as the state space of the random walk.
Only results are presented here since the detailed proofs will appear
elsewhere.
\end{abstract}

\tableofcontents
\section{Introduction}
\label{sec_intro}
\subsection{On the generality of the random walk in a random environment
on a tree}
\label{ssec_genericity}
Markov chains on denumerable graphs enter in the modelling
of a rich variety of phenomena; among such graphs, trees
play a basic and generic r\^ole in the sense that they can encode
simultaneously
\begin{itemize}
\item
 the topological structure of a vast class of general (\textit{i.e.}
not necessarily tree-like) denumerable
graphs,
\item 
the combinatorial structure of paths on these general graphs, and
\item
the probability structure generated on the trajectory space
of the Markov chain evolving on the these general graphs.
\end{itemize} 
In the sequel we give a short explanation of the reason trees play such a 
generic rôle among denumerable graphs. We know, after Kolmogorov, 
that in order to describe a randov variable $X$ with space of outcomes
a discrete measurable set $(\BbX,\cX)$, of law $\BbP_X$, one has
to use an abstract probability space $(\Omega,\cF,\BbP)$ on which the
randov variable is defined. However, the choice of this space is not unique.
Among the infinite possible choices, there exists a ``minimal'' one,
given by $\Omega=\BbX$, realising the random variable as the identity map
$X(\omega)=\omega$ and $\BbP=\BbP_X$. Similarly, when $(X_n)_{n\in\BbN}$
is a sequence of \textit{independent} and identically distributed 
random variables with space of outcomes  the discrete space $(\BbX, \cX)$,
the minimal realisation of the abstract probability space carrying the
whole sequence is  the \textit{space of trajectories or full
shift} $\Omega=
\BbX^\BbN$ and the sequence is realised by the canonical projection
$X_n(\omega)=\omega_n$, for $n\in\BbN$. However for 
a sequence of \textit{dependent} random variables the space of trajectories
may be not minimal. The reader can easily  convince herself by considering
the example of $\BbX=\{a,b,\ldots, z\}$ and  $(X_n)_{n\in\BbN}$ being 
the sequence
of letters appearing in a natural language. Then the occurance 
$(X_n, X_{n+1}, X_{n+2})=rzt$ never appears in a language like English or
French. 

For the special case of Markovian dependence, the natural space
for the realisation of the sequence $(X_n)_{n\in\BbN})$ is the so called
unilateral or bilateral \textit{subshift spaces} obtained from thh full
shift by deleting all sequences containing forbidden subwords. For Markovian
sequences, defined through a stochastic matrix $P$, 
subshift spaces can also be obtained in terms of the
adjacency matrix of the sequence $A:\BbX\times\BbX\rightarrow \{0, 1\}$,
given by the formula $A(x,y)=1$ whenever $P(x,y)>0$ and zero otherwise.

Suppose now we are given an arbitrary directed graph
having  a finite or denumerable set of vertices and being such that
only a finite number of edges is emitted out of every vertex, not having 
any sources or sinks, and
having (complex) weights attached on its edges.
Since the vertex set is at most countable, it is in bijection with an at
most countable alphabet. We say that a sequence of letters from the
alphabet is a path of the graph, if the corresponding sequence
of vertices is such that any ordered pair of subsequent vertices
is an element of the edge set.  The set of paths of arbitrary length is 
called \textit{path space
of the graph} (see
\cite{KumPasRaeRen} for precise definitions.) Notice that the condition
that two subsequent vertices must be an allowed edge prevents some
sequences of vertices from being a path.
The adjacency matrix defines a subshift space for the
trajectories of a symbolic unilateral or bilateral 
Markov chain that is identifed with the  path space
of the graph.
It is easy to show that
the path space of any graph
has a natural tree structure
\cite{CamPet-wien}, giving thus a first hint that trees play a prominent
rôle among graphs. Their importance does not stop here however. We can in fact
define  a so-called
\textit{evaluation map} from the path space of a graph
into some set $\BbA$. Depending on the precise algebraic structure
of the set $\BbA$ and of the evaluation map, a vast class of objects can be 
defined. Instead of giving precise definitions, let us give
the following
\begin{exam}{(An elementary case)}
Start with the finite complete graph on 4 vertices, denote them 
by the letters E, N, W, and 
S for definiteness, and consider the unilateral 
path space of the graph, \textit{i.e.}
the set of words of arbitrary length on the alphabet of these 4 letters.
Choose for the space  $\BbA$ the Abelian group $\BbZ^2$ and for the evaluation
map the function taking the value
$n_E \be_1+n_N \be_2+ n_W (-\be_1)+ n_S (-\be_2)\in \BbZ^2$, where
$n_E,n_N, n_W, n_S$ denote the occurrences of the letters 
E, N, W, and S in a given word, and $\be_1,\be_2$ are the unit vectors
of $\BbZ^2$. 
The evaluation map can be thought as the physical obervable ``position of the
random walker'' when the history of the individual directions of movement 
is kept
in memory; individual directions 
define a path in the path space and the evaluation map
computes the actual position of the walker which has done
$n_E$ eastward movements, $n_W$ westward movements, etc.   
The graph that is generated in this way is
the graph that coincides with
the Cayley graph of the 
Abelian group $\BbZ^2$. The evaluation map being many-to-one,
its multiplicity encodes the combinatorics of the path space
while the product of the graph weights gives the relative
weight of individual paths in the path space of this
graph. When the weights are probability
vectors, this weight coincides with the probability of the
trajectory of the simple random walk on  $\BbZ^2$. Notice however
that the weights need not to be probability vectors, allowing us
to consider both random and quantum grammars.
\end{exam}
The above example serves as a basic paradigm; by appropriately
changing the evaluation map, it can be generalised in numerous ways 
to produce configuration space for DNA strands (when the path space
coincides with the path space of the complete graph over
4 letters), non-Abelian groups like the free group on 2 generators, 
Cuntz-Krieger
$C^*$-algebras (\textit{i.e}\ non-commutative operator algebras
attaching non-zero partial isometries on every edge; 
see \cite{Pet-GDG} for additional details) naturally 
arising in some problems of
quantum information.
Infinite trees provide
a rich variety of mathematical problems, particularly connected to
their non-amenability; beyond their mathematical interest they arise
as more or less realistic models in several applied fields like random
search algorithms in large data structures, Internet traffic, random
grammars and probabilistic Turing machines, DNA coding, interacting
random strings and automated languages, etc.
The previously exposed ideas, already lurking in \cite{Mal-ISS,Mal-SEGG}, have
been exemplified in \cite{FlaSed,CamPet-wien}, and will be further
exploited in \cite{Pet-GDG}. 

Random environment is a means to introduce context-dependence in the language.
Again, due to the generality of trees as underlying graphs, random
walks in random environment on trees provide interesting context-dependent
results for random walks on a huge class of more general graphs.

Random walks
in random environments on various types of graphs are
known to display a behaviour dramatically differing from the one for
ordinary random walks on the same graph.
See \cite{Sol,Sin,KesKozSpi,Kes,LyoPem,ComMenPop,Mal-RG,MenPet-rwre,ComPop,CamPet-rwrol} etc.\
for a very partial list of known results and models.

\subsection{Rough statement of the main results for bindweeds}

As a simple example, consider a rooted tree with constant branching $b$,
and let us construct a random walk in random environment on it by
sampling the transition probabilities (or transition rates, for
continuous time) in each vertex from a given distribution,
independently. When studying the question of (positive) recurrence
of such a walk, one naturally arrives (see~\cite{LyoPem,MenPet-rwre})
at the following model. On each edge $a$ of the tree, we place
an independent copy of a
positive random variable ${\hat \xi}_a$. Then, for each vertex $\bv$
denote ${\hat \xi}[\bv]={\hat \xi}_{a_1}\ldots{\hat \xi}_{a_n}$,
where $a_1,\ldots,a_n$ is the (unique) path connecting the root to $\bv$.
Models of this type are called {\it multiplicative cascades\/}
(see e.g.~\cite{KahPey,LiuRou}), and, to study the positive recurrence
of the corresponding random walk in random environment, one has to
answer the question whether the sum of~${\hat \xi}[\bv]$ is finite.
It turns out (see e.g. \cite{LyoPem}) that the classification
parameter for this problem is
\begin{equation}
\label{class_param}
\hat\lambda = \inf_{s\in[0,1]}\BbE {\hat \xi}_a^s,
\end{equation}
which is then compared to $1/b$ (in fact, in~\cite{LyoPem} the case
of general tree was considered) and the following result is established:
\begin{itemize}
\item
if $\hat\lambda b<1$ then the random walk for almost all environments is positive
recurrent
\item
if $\hat\lambda b>1$ then the random walk for almost all environments is 
transient.
\end{itemize}
The critical case $\hat\lambda b=1$ is more complicated. 

In \cite{MenPetPop}, we introduced a model of random walk  with internal 
degrees of freedom whose multiplicative cascade counterpart is
expressed by  a model placing  
random matrices on the
edges in place of scalar
random variables. The main classification parameter, $\lambda$, will be
defined in the formula (\ref{eq-k(s)}) below, and the main results are
Theorems \ref{th-mmc-ergodicity} and \ref{th-mmc-transience}
(a lot of preliminary work is required, however, before formulating
these results).   This model
is also equivalent to  random walk in random
environment on a multiplexed tree, a process we term \textit{bindweed}
in the sequel. However, for trees with average branching $b$,
the rough statement of our result is as above with $\hat\lambda$ replaced
by $\lambda$.

It was remarked in \cite{MenPet-rwre} that asymptotic properties
like recurrence/transience of random walk on trees with constant
branching $b$ are intimately connected to the existence of non-trivial
solutions for the so-called multiplicative chaos equation of order~$b$,
first introduced as a simple turbulence model in \cite{Man}.
The simplest variant of the multiplicative chaos equation is the following:
let $(\xi_i)_{i=1,\ldots, b}$, with $b\in\BbN$, be a finite family
of non-negative random variables having known joint distribution and
$(Y'_i)_{i=1,\ldots, b}$ and $Y$ be a family of $b+1$ independent non-negative
random variables distributed according the same unknown law and verifying
\[
Y\eqlaw\sum_{i=1}^b Y'_i \xi_i.
\]
The multiplicative chaos problem consists in determining under which
conditions on the joint distribution of the $\xi$'s the above equation has
a non-trivial solution.
This scalar problem is thoroughly studied in the
literature, see e.g.\ \cite{ColKou,DurLig,Liu97}.
As we remark later in this paper, the matrix multiplicative chaos equation
may be an interesting problem to study as well.

\section{Notation}
\label{ssec-notations}
In this section we give the formal definitions concerning trees, in particular,
we define the notions of the growth rate and the branching number.

We denote $\BbR_+=[0,\infty[$, $\PosInt=\{0,1,2,\ldots\}$,
$\StrPosInt=\{1,2,3,\ldots\}$, and
for every $n\in \StrPosInt$, $\BbN_n=\{1,2,\ldots,n\}$,
while $\BbN_0=\emptyset$.
Let $\cA\equiv\cA^1$ be a finite or infinite denumerable set, called
the \textit{alphabet}. Define $\cA^0=\{\emptyset\}$ and
for every $n\in \StrPosInt$ denote
\[
\cA^n=\{\bom{\alpha}=\alpha_1\cdots\alpha_n:
\alpha_i\in \cA\ \textrm{for}\ i\in\BbN_n\}
\]
the set of words of length $n$ (\textit{i.e.}\ having $n$ letters),
\[
\cA^*=\cup_{n\in\PosInt}\cA^n
\]
the set of words of arbitrary (finite) length,
and
\[
\partial\cA^*\equiv\cA^\infty=
\{\bom{\alpha}=\alpha_1\alpha_2\cdots: \alpha_i\in \cA\ \textrm{for}\
i\in\StrPosInt\}
\]
the set of infinite words.
Finally, denote $\ol{\cA^*}=\cA^*\cup\partial\cA^*$ and
$\oc{\cA^*}=\cA^*\setminus \cA^0$.

For every $\bom{\alpha}\in\cA^*$, there exists $n\in\PosInt$ such that
$\bom{\alpha}\in \cA^n$; in this situation $|\bom{\alpha}|:=n$ denotes
the \textit{length} of the word $\bom{\alpha}$ with the convention
$|\emptyset|=0$. Consistently, for every $\bom{\alpha}\in\partial\cA^*$,
we have $|\bom{\alpha}|=\infty$. For $\bom{\alpha}\in\ol{\cA^*}$ with
$|\bom{\alpha}|\geq n$ we denote by
 $\bom{\alpha}\rt{n}=
\alpha_1\cdots\alpha_n\in\cA^n$ the
\textit{restriction} of $\bom{\alpha}$ to its $n$ first letters
with the convention $\bom{\alpha}\rt{0}=\emptyset$.
For every $\bom{\alpha}\in\oc{\cA^*}$, the \textit{ancestor}
$\hat{\bom{\alpha}}$ of $\bom{\alpha}$ is defined
by $\hat{\bom{\alpha}}=\bom{\alpha}\rt{|\bom{\alpha}|-1}$.
For $\bom{\alpha}\in\cA^*$ and $\bom{\beta}\in\ol{\cA^*}$, the
\textit{concatenation} of $\bom{\alpha}$ followed by
$\bom{\beta}$ is the word $\bom{\alpha}\bom{\beta}=
\alpha_1\cdots\alpha_{|\bom{\alpha}|}\beta_1\beta_2\cdots$
and for $\bom{\alpha},\bom{\beta}\in \cA^*$, their \textit{common radix}
$\bom{\alpha}\wedge\bom{\beta}$ is the longest word
$\bom{\gamma}\in\cA^*$ such that $\bom{\alpha}=\bom{\gamma}\bom{\alpha'}$
and $\bom{\beta}=\bom{\gamma}\bom{\beta'}$ for some words \
$\bom{\alpha'},\bom{\beta'}\in \cA^*$. We write
$\bom{\alpha}\leq\bom{\beta}$ if
$\bom{\alpha}=\bom{\alpha}\wedge\bom{\beta}$.

\Rk Notice that, consistently with the above notation, the symbol
$\BbN^*$ denotes the set of finite words on the alphabet $\BbN$, contrary
to some tradition (especially the French one)
where this symbol is used to denote what we call here $\StrPosInt$.

\begin{defi}{}
A mapping $B:\BbN^*\rightarrow\PosInt$ is
called a \textit{branching function}.
\end{defi}

To each branching function corresponds a uniquely
determined  rooted tree $\BbT=(\BbV,\BbA)$ with
vertex set $\BbV\equiv\BbV^*(B)\subseteq\BbN^*$ and edge set
$\BbA=\oc{\BbV}$ defined as follows:
$\BbV^*(B)=\cup_{n\in\BbN}\BbV^n(B)$ where $\BbV^0(B)=\{\emptyset\}=
\{\textrm{root}\}\equiv\{\0\}$ and for $n\in\StrPosInt$,
\[
\BbV^n(B)=\{\bv=v_1\cdots v_n: v_l\in \BbN_{B(\bv\rt{l-1})},
\ \textrm{for}\ l=1,\ldots, n\}.
\]
The branching function is said to be \textit{without extinction}
if the corresponding
tree has non-trivial boundary $\partial\BbV$.
The edge set is the subset of unordered pairs of vertices
$[\bu,\bv]=[\bv,\bu]$ such that either $\bv=\hat{\bu}$ or
$\bu=\hat{\bv}$. Since every vertex has a unique ancestor, every edge
is indexed by its outmost vertex, \textit{i.e.}\ for
every $\bv\in\oc{\BbV}$, the corresponding edge is $a(\bv)=[\hat{\bv},\bv]$,
showing thus that $\BbA\simeq \oc{\BbV}$.

If $\bu,\bv\in\BbV$ and $\bu\leq\bv$ we define the
\textit{path} $[\bu,\bv]$ as the  collection of the $|\bv|-|\bu|$ edges
$[\bu, \bv\rt{|\bu|+1}],\ldots, [\hat{\bv},\bv]$, and if $\bu=\emptyset$
then we simply denote by $[\bv]$ the path $[\emptyset,\bv]$ for every
$\bv\in\oc{\BbV}$.
In the sequel we shall consider only \textit{branching
functions without extinction}.
\begin{defi}{}
\label{def-grr}
Let $\kappa_n=\card \BbV^n(B)$
denote the cardinality of the $n^{\textrm{th}}$
generation of the tree defined by the branching function without
extinction $B$.
We call \textit{lower growth rate} of the tree
\[
\lgr=\liminf_n \kappa_n^{1/n},
\]
\textit{upper growth rate} of the tree
\[
\ugr=\limsup_n \kappa_n^{1/n},
\]
and, if $\lgr=\ugr$, we call the common value \textit{growth rate}
\[
\grr=\lim_n \kappa_n^{1/n}.
\]
\end{defi}
For $\bu,\bv\in\partial\BbV$, define $\delta(\bu,\bv)=\exp(-|\bu\wedge\bv|)$.
It can be shown that $\delta$ is a distance on $\partial\BbV$. Moreover
if $\|B\|_\infty=\sup_{\bv\in\BbV} B(\bv)<\infty$ then the space
$(\partial\BbV,\delta)$ is compact and
we can define its Hausdorff dimension
$\dim_H\partial\BbV$ as usual (see \cite{Falconer} for instance).
\begin{defi}{}
\label{def-brr}
For a tree $\BbV$ generated by a branching function $B$ with
$\|B\|_\infty<\infty$, we define its \textit{branching rate}
\[
\brr=\exp(\dim_H\partial\BbV).
\]
\end{defi}
It is shown in \cite{GraMauWil,LyoPem} that
$\brr=\sup\{\lambda: \inf\sum_{v\in C}\lambda^{-|\bv|}>0\}$
where the infimum is evaluated over all cutsets $C$ of $\BbV$. We have
in general that $\brr\leq \lgr$.

\section{Matrix multiplicative cascades and the corresponding results}
\label{ssec-mmc}
Let $(\Omega,\cF,\BbP)$ be some abstract probability space which carries
all the random variables that will be needed in the model. Let $(\BbV,\BbA)$
be the rooted tree associated with a given branching function $B$. Let
$G$ be the topological group $\GLdR$, $\cG$ its Borel $\sigma$-algebra
and $\mu$ a probability on $(G,\cG)$. Denote by
$\sigma_\mu=\supp \mu\subset G$ the support
of the measure $\mu$ and by $\Sigma_\mu$ the semi-group
generated by $\sigma_\mu$. On $(\Omega,\cF,\BbP)$, define
an edge-indexed family of independent $G$-valued random variables
$(\xi_a)_{a\in\BbA}$ identically
distributed according to $\mu$, \textit{i.e.}\
\[
\BbP(\xi_a\in \dd g)=\mu(\dd g), \mbox{ for all } a\in\BbA.
\]
For $\bu, \bv\in\BbV$ with $\bu\leq\bv$ define
\[
\xi[\bu,\bv]\equiv\prod\limits^{\leftarrow}_{a\in[\bu,\bv]}\xi_a,
\]
where $\prod\limits^{\leftarrow}$ denotes the product in reverse order,
\textit{i.e.}\ if $[\bu,\bv]=a_1\cdots a_k$ then
$\xi[\bu,\bv]=\xi_{a_k}\cdots \xi_{a_1}$ with the
convention $\xi[\bv,\bv]=e$ where $e$ is the neutral element
of $G$. We introduce the following
$G$-valued random processes:
the \textit{matrix-multiplicative cascade process}
\[
\psi_n=\sum_{\bv\in\BbV^n} \xi[\bv], \quad n \in \BbN
\]
and the \textit{integrated matrix-multiplicative cascade process}
\[
\zeta_n=\sum_{k=1}^n \psi_k, \quad n \in\StrPosInt.
\]

For a fixed $\bv\in\partial \BbV$ and all $n\in\BbN$
\[
X_n\equiv X_n(\bv)= \xi[\bv\rt{n}].
\]
It is immediate to see (cf.~\cite{Revuz})
that $(X_n)_{n\in\StrPosInt}$ is a $G$-valued
multiplicative Markov chain with stochastic kernel
\[
P(g,\dd g')\equiv\BbP(X_{n+1}\in \dd g'|X_n=g)=\mu\star\delta_g(\dd g'),
\quad n\in \StrPosInt,
\]
where for two measures $\mu,\mu'$ on $(G,\cG)$ their convolution
$\mu\star\mu'$ is defined by its dual action on $L^1(G)$ via
\[
\bra \mu\star\mu', f\ket \equiv \int_G f(g) \mu\star\mu'(\dd g)=
\int_G\int_G f(gg')\mu(\dd g)\mu'(\dd g'),
\]
for all $f\in L^1(G)$.
We equip $\End(\BbR^d)\equiv \GLdR=G$ with the operator
norm, denoted $\|\cdot\|$, stemming from the $l_1$ norm of the vector
space $\BbR^d$.

In order to be able to  apply the results of \cite{GuiLeP} to our special case,
we require the following conditions on $\mu$:

\noindent
\textbf{Condition 1 (Integrability):}\
For all $s\in\BbR_+$,
\[
 \int_G \|g\|^s \mu(\dd g)<\infty.
\]

\noindent
\textbf{Condition 2 (Strong irreducibility):}\
 We assume that the set $\Sigma_\mu$
is strongly irreducible, \textit{i.e.}\ there is no finite 
$\Sigma_\mu$-invariant family of
proper   subspaces.

\noindent
\textbf{Condition 3 (Strict positivity):}\
We assume that
$\sigma_\mu\subseteq \ol{G}_+$
and that
$\mu(\ol{G}_+\setminus G_+)=0$.

For $s\geq 0$, define
\begin{equation}
\label{def-k(s)}
k(s)=\lim_n\left(\int_G \|g\|^s\mu^{\star n}(\dd g)\right)^{1/n}.
\end{equation}
(By virtue of theorem 1 of \cite{GuiLeP}, 
under conditions 1--3 this limit exists in $\BbR_+$ and defines
a log-convex function. As a matter of fact, in \cite{GuiLeP} a
weaker condition than 3, called proximality, is needed to prove this result.)
We define in the sequel the quantity $\lambda$, that turns out to be
the main classification parameter for the matrix multiplicative
cascades model, by
\begin{equation}
\label{eq-k(s)}
\lambda= \inf_{s\in[0,1]} k(s)
\end{equation}
(compare (\ref{eq-k(s)}) with (\ref{class_param}).)

We are now in the position to state our main results.

\begin{theo}
\label{th-mmc-ergodicity}
Let $(\BbV,\BbA)$ be some tree defined in terms of
a given  branching function $B$ and $\ugr$ and $\lambda$
defined as in definition \ref{def-grr} and
equation (\ref{eq-k(s)}) respectively.
Under the conditions 1, 2, and 3,
\[
\lambda \ugr <1 \Rightarrow
\zeta_{\infty, ij}<\infty
 \textrm{ almost surely, for all}\ i,j=1,\ldots,d.
\]
\end{theo}

\begin{theo}
\label{th-mmc-transience}
Let $\brr$ and $\lambda$ be the quantities introduced
in definition \ref{def-brr} and equation \ref{eq-k(s)}
respectively and let
$\chi\in \BbR^d$ be the vector having all its components equal to 1:
$\chi_i=1$, for all $i=1,\ldots, d$.
Let $(\BbV,\BbA)$ be some tree defined in terms of
a given branching function~$B$ without extinction.
Under the conditions 1, 2, and 3',
\[
\lambda \brr >1\Rightarrow Z_\infty := (\chi, \zeta_\infty\chi)=
\infty \ \textrm{almost surely.}
\]
\end{theo}

\Rk Similarly to \cite{LyoPem}, there is a gap between
Theorems \ref{th-mmc-ergodicity} and \ref{th-mmc-transience}, since in
general the branching number need not be equal to the growth rate.
However, this is not very important, because in most of the practical
examples these quantities do coincide.

\Rk As mentioned above, the classification parameter for this problem is
$\lambda=\inf_{s\in[0,1]} k(s)$. This parameter is not
explicitly computable in general since it involves the infinite product of
matrices. However for some particular cases this quantity
can be computed explicitly as stated in the following proposition.

\begin{prop}
\label{prop-lambda=rho}
 Suppose that the measure
$\mu$ is such that $g_{ij}<1/d$ almost surely for all $i,j=1,\ldots,d$.
Then $\lambda$ is the largest eigenvalue of the matrix $\BbE g$.
\end{prop}

\Rk It is interesting to consider the \textit{chaos equation}
for the case of matrix-valued random variables and constant branching $b$:
\begin{equation}
\label{chaos_eq}
Y \eqlaw \sum_{j=1}^b  Y'_j\xi_j,
\end{equation}
where $Y,Y'_j,\xi_j$ are $G$-valued random variables, and $\xi_j$
(which are not necessarily independent) are
distributed according to~$\mu$;
$Y'_j$, $j=1,\ldots,b$, are i.i.d.\ and
have the same (unknown) law as $Y$.
Analogously to \cite{MenPet-rwre} we
can get that (at least in the case when $\xi_1$ satisfies
conditions 1,2,3)
$\lambda d = 1$ is a necessary condition for the
existence of solution of (\ref{chaos_eq}). It is an open
problem whether this condition is sufficient.

\Rk The condition of independence of the random variables $\xi_a$
can be relaxed; what is important is 
\begin{enumerate}
\item
if $\xi_a$ and $\xi_b$
are not adjacent to the same vertex then they must be independent,
and 
\item
the $\xi$'s that belong to any path emanating from the
root must be independent.
\end{enumerate}

\section{The bindweed model}
\label{ssec-rw}
In this section we introduce a model describing an evolution
of a random string in random environment on a tree
(which is somewhat similar to the model studied in \cite{ComMenPop})
which we call the bindweed model. Then, we show that its
classification from the point of view of positive recurrence
can be obtained by using theorems \ref{th-mmc-ergodicity}
and \ref{th-mmc-transience}.

Let $\cS
=\{1,\ldots,d\}$ be a finite alphabet and denote,
in accordance with the notations introduced in Section~\ref{ssec-notations},
$\cS^{n+1}=\{\bom{\sigma}=\sigma_0 \cdots \sigma_n: \sigma_i\in\cS\}$
the
set of words of length $n+1$ composed from the symbols of the alphabet $\cS$,
$\cS^0$ the set containing only the empty word and $\cS^*$
the set of words of arbitrary length. Suppose that a branching function $B$
is given on $\BbN^*$ and denote $\BbV^n\equiv\BbV^n(B)$ the corresponding
generations of the tree determined by $B$. Therefore, the
rooted tree $\BbT=(\BbV,\BbA)$
is uniquely defined.

\begin{figure}
\centerline{\includegraphics[width=0.5\textwidth]{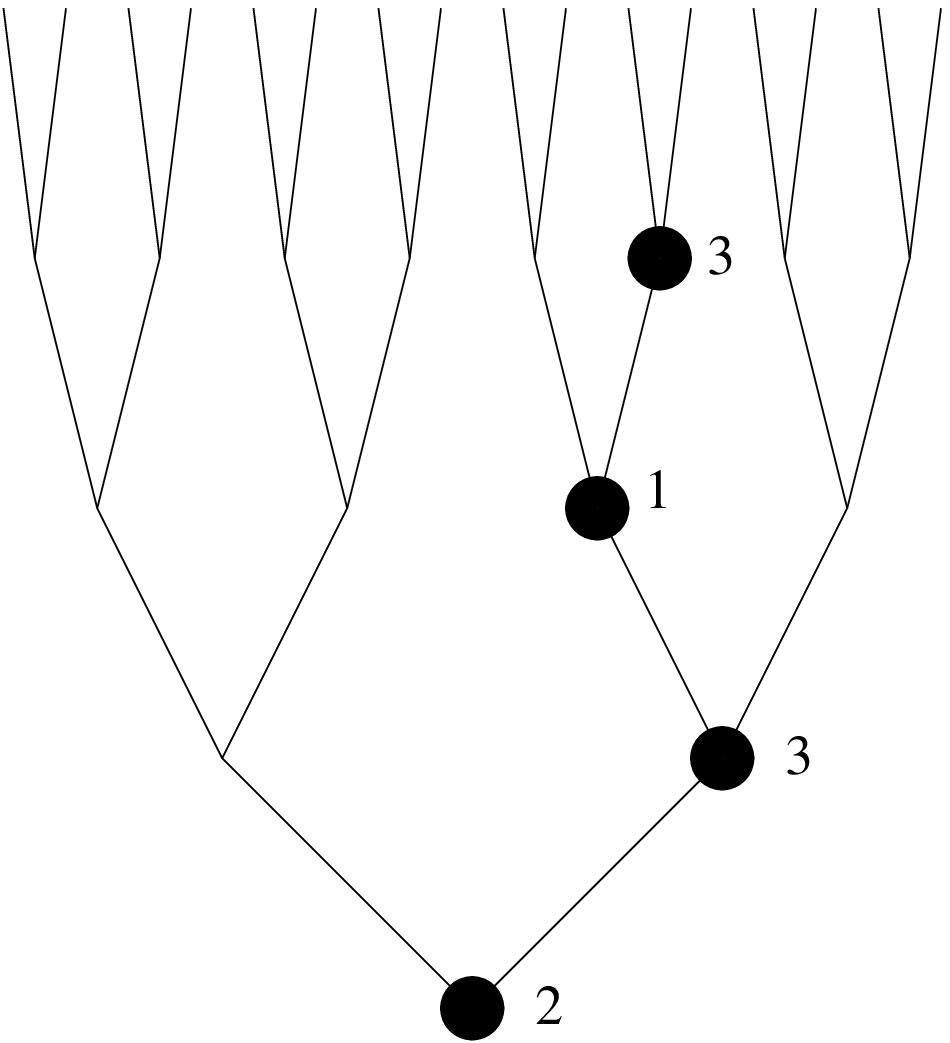}}
\caption{\textsf{A typical state of the bindweed model for $\cS=\{1,2,3\}$,
$\bom{\sigma}=2313$, and $\BbT$ the  binary tree.}}
\label{fig_bindweed}
\end{figure}

Now
we are going to construct a continuous-time Markov chain
with state space $\mathfrak{S}$, defined by
\[
\mathfrak{S} = \{\hat{\emptyset}\}\cup \bigcup_{n=1}
           (\BbV^n\times\cS^{n+1}),
\]
where $\hat{\emptyset}$ is a special state to be defined
later. In fact, what happens is the following: we place
a word $\bom{\sigma}=\sigma_0\ldots \sigma_n$ on the tree $\BbT$
in such a way that the $0^{\textrm{th}}$ symbol of the word is placed
on the root $\0$, for any $i=1,\ldots,n$ the $i^{\textrm{th}}$ symbol
of the word is placed somewhere in $\BbV^i$,
and, if the $i^{\textrm{th}}$ symbol $\sigma_i$ is placed
on vertex $\bu$, and $\sigma_{i+1}$ on $\bv$, then
$\bu<\bv$ and $[\bu,\bv]\in\BbA$ (see figure \ref{fig_bindweed}).
The state ${\hat \emptyset}$ means that nothing
is placed on the tree.

Now, let us define the dynamics of the bindweed model.
Suppose that for any $a\in\BbA$ two collections of positive numbers
$(\nu_{yz}(a), y,z \in \cS )$, $(\mu_y(a), y \in \cS )$ are given.
If the bindweed model is in the state $(\bu,\bom{\sigma})$, where
$\bom{\sigma}=\sigma_0\ldots \sigma_{n-1}y$, $\bu\in\BbV^n$, then
\begin{itemize}
\item for $n\geq 0$ it jumps to the state
$(\bv, \sigma_0\ldots \sigma_{n-1}y z)$
 with rate $\nu_{y z}(a(\bv))$, for all $\bv\in\BbV: \bu=\hat{\bv}$;
\item for $n\geq 1$ it jumps to the state $(\hat{\bu}, \sigma_0\ldots
\sigma_{n-1})$
with rate $\mu_y(a(\bu))$.
\end{itemize}
For any $\sigma_0\in\cS$ the transitions
${\hat \emptyset}\to (\0,\sigma_0)$
and $(\0,\sigma_0)\to {\hat \emptyset}$ occur with rate $1$.
Thus, we have defined a continuous-time Markov chain
with  state space $\mathfrak{S}$.

Let us describe now how to choose the transition rates.
Let $\rho$ be any probability measure on $\BbR^{d^2+d}_+$.
Suppose that for any $a\in\BbA$ the vector
$\Xi(a)=(\nu_{yz}(a),y,z\in\cS,
\mu_y,y\in\cS)$ is random, having distribution $\rho$,
and $(\Xi(a), a\in\BbA)$ are independent and identically distributed.
Fix a realisation
of that collection of random vectors and consider
the bindweed model with the transition rates ruled by that
realisation.
So, the model that we constructed is a continuous-time Markov
chain in a quenched random environment.

Now, we are interested in obtaining a classification of
this Markov chain with respect to positive recurrence.
For $(\bv,\bom{\sigma})\in\mathfrak{S}$ denote by
$\pi(\bv,\bom{\sigma})$ the stationary measure.
For any $a\in\BbA$ let $\xi_a$ be a $d\times d$ matrix
whose matrix elements are
defined in the following way: $\xi_{a, xy}=\nu_{xy}(a)/\mu_y(a)$,
$x,y\in\cS$. It is not difficult to see that
 we have a reversible Markov chain, so
it is clear that $\pi({\hat \emptyset}) = \pi(\0,x)$,
for all $x\in\cS$, and, for any $\bv\in\BbV^n$, $n\geq 1$,
and $x,y,\sigma_0,\ldots,\sigma_{n-2}\in\cS$,
we can formally write
\begin{eqnarray}
\pi(\bv,\sigma_0\ldots \sigma_{n-2}x y)
& = & \frac{\nu_{xy}(a(\bv))}{\mu_y(a(\bv))}
           \pi(\hat{\bv},\sigma_0\ldots \sigma_{n-2}x)\nonumber \\
        &=&  \xi_{a(\bv),xy} \pi(\hat{\bv},\sigma_0\ldots \sigma_{n-2}x)
\label{eq-revers}.
\end{eqnarray}

Then it is shown in \cite{MenPetPop} that 
\[
\sum_{\substack{\bv\in\BbV^n\\\bom{\sigma}\in \cS^n}}
\pi(\bv,\bom{\sigma}) = \pi({\hat \emptyset})
 \left(\chi, \sum_{\bv\in\BbV^n} \xi[\bv] \chi\right)
\] where
$\chi$ is the  vector of order $d$ with all its coordinates
equal to $1$.
Thus $\sum_{(\bv,\bom{\sigma})
\in\mathfrak{S}}\pi(\bv,\bom{\sigma})$
is finite if and only if $Z_{\infty}$
is finite. Thus, theorems \ref{th-mmc-ergodicity} and \ref{th-mmc-transience}
allow us to obtain the
classification of the bindweed model in random environment
from the point of view of positive recurrence, in the following way:
\begin{prop}
Suppose that the distribution of the random matrix $\xi_a$ is
such that the Conditions 1, 2, and 3 are satisfied.
Let $\lambda$ be the quantity defined as in Section~\ref{ssec-mmc}.
Then
\begin{itemize}
\item if $\lambda\ugr < 1$, then the bindweed model is positive recurrent;
\item if $\lambda\brr > 1$, then the bindweed model is not positive recurrent.
\end{itemize}
\end{prop}

\section{Open problems and further developments}
We demonstrated a close relationship between matrix multiplicative
cascades and random walks in random environment on multiplexed trees.
In particular it is proven in \cite{MenPetPop} that both systems are classified
by the same parameter.
However, the critical region remains out of reach for the moment. Firstly
it is not known whether, for sufficiently recular trees so that
$\ugr=\brr=b$, the walk is null recurrent for $\lambda b=1$ or
some additional condition is needed on $\mu$ 
as is the case for scalar
multiplicative chaos \cite{Liu97} and for random walk \cite{LyoPem}.
Returning to the general tree where  $\ugr\not=\brr$, we obtain a gap
in the space of classifying parameters. It is however conjectured in 
\cite{Mal-ISS} that in general, the set of critical values is of zero Lebesgue 
measure for string problems. We expect the same phenomenon
to occur here. Nevertheless, whether
the critical value of the  parametre is $\brr$ or $\ugr$ or some intermediate
value is unknown for the moment.

An important step towards understanding these problem should be made
if conditions for the existence of non-trivial fixed
points of the functional equation (\ref{chaos_eq}) were obtained. This
remains for the moment an open problem although under investigation.

\bibliographystyle{plain}
\bibliography{rwre,matrix,petritis,cstar}
\end{document}